\documentclass[12pt]{article}
\usepackage{amssymb,amsfonts,amsmath}
\usepackage{graphics}
\usepackage{graphicx}
\usepackage{epsfig}
\usepackage{color}
\usepackage{enumerate}

\newcommand\CC{{\mathbb C}}
\newcommand\GA{\mathcal A}

\newcommand\LL{{\cal L}}

\newcommand\NN{{\mathbb N}}
\newcommand\RR{{\mathbb R}}

\newcommand\GX{{\mathcal X}}

\newcommand\Id{{I}}

\def\beq{\begin{equation}}
\def\eeq{\end{equation}}
\newtheorem{thm}{Theorem}[section]

\newtheorem{lem}[thm]{Lemma}
\newtheorem{cor}[thm]{Corollary}
\newtheorem{rem}[thm]{Remark}

\newtheorem{ex}[thm]{Example}

\newcommand\re{\mathop{\rm Re}\nolimits}

\newcommand\supp{\mathop{\rm supp}\nolimits}

\newcommand\rad{\mathop{\rm Rad}\nolimits}

\def\C{{\mathbb C}}

\def\R{{\mathbb R}}

\newcommand{\ds}{\displaystyle}

\newcommand{\g}{\gamma}

\renewcommand\phi{\varphi}
\newcommand\FB{\mathcal{F}\mathcal{B}}

\title{Dichotomy results for norm estimates in operator semigroups}
\author{\rm I. Chalendar\thanks{Universit\'e de Lyon; CNRS; Universit\'e Lyon 1; INSA de Lyon; Ecole Centrale de Lyon,  CNRS, UMR 5208, Institut Camille Jordan, 43 bld. du 11 novembre 1918, F-69622 Villeurbanne Cedex, France. \tt  chalendar@math.univ-lyon1.fr},
 J. Esterle\thanks{I.M.B., UMR 5251, Universit\'e de Bordeaux,
351
cours de la Lib\'eration,
33405 Talence Cedex, France.
  \protect\linebreak[3]
{\tt jean.esterle@math.u-bordeaux1.fr}.}
\ and   J.R. Partington\thanks{School of
Mathematics,
University of Leeds, Leeds LS2 9JT, U.K.
\protect\linebreak[3]
{\tt J.R.Partington@leeds.ac.uk}.}
}

\date{}
\begin{document}

\baselineskip18pt

\maketitle

\bibliographystyle{plain}

\begin{center}
{\it  Dedicated to Charles Batty on the occasion of his sixtieth birthday.}\\
\end{center}

\begin{center} {\it
In line with an Oxford tradition,\\
Charles Batty was sent on a mission:\\
to instruct all the troops\\
in $C_0$ semigroups –-\\
this indeed was a brilliant decision!
}
\end{center}

\section{Introduction}

We recall that a one-parameter family $(T(t))_{0 <
  t<\infty}$   in a Banach algebra (often itself simply the algebra of bounded linear operators on a Banach space $\GX$)
is a semigroup if 
\[T(s+t)=T(s)T(t)\qquad \mbox{ for all }t,s> 0.
\]
We shall be concerned here with semigroups that are strongly continuous on $\RR_+:=(0,\infty)$, but not necessarily norm-continuous at the origin. 
As an example to bear in mind we mention the semigroup $T(t):x \mapsto x^t$ in the algebra $C_0([0,1])$ of continuous functions on $[0,1]$
vanishing at the origin, which will be discussed later.

Later, we consider semigroups defined on a sector in the complex plane, in which case they will be assumed to be {\em analytic\/}: that is, complex-differentiable in the norm
topology.
\\

The results in this survey indicate that the quantitative behaviour of the semigroup at the origin provides
additional qualitative information, such as uniform continuity or analyticity. Here are a few examples.\\

We recall the classical {\em zero--one law}, asserting that
if $\limsup_{t \to 0^+} \|T(t)-\Id\|<1$, then in fact $\|T(t)-\Id\| \to 0$ and hence the semigroup is
uniformly continuous, and of the form $e^{tA}$ for some bounded operator $A$.
To see this,
set $L=\limsup_{t\to 0^+ } \|T(t)-\Id\|$. 
 Since
$$2(T(t)-\Id)=T(2t)-\Id-(T(t)-\Id)^2,$$ 
we have $2L\leq L+L^2$,
 and thus $L=0$ or $L\geq 1$. \\

Another result involving the asymptotic behaviour at $0$ and providing a uniformly continuous semigroup is the 
following, proved in 1950 by Hille~\cite{hille50}
(see also \cite[Thm.~10.3.6]{HP74}).  This result is usually stated for $n=1,$ but Hille's argument works for any positive integer.

\begin{thm}\label{thm:hille50} Let $(T(t))_{t>0}$ be a $n$-times continuously differentiable semigroup over the positive reals. If $\limsup_{t\to 0^+}\Vert t^nT^{(n)}(t)\Vert <({n\over e})^n,$ then the generator of the semigroup is bounded.
\end{thm}

In the direction of analyticity, a classical result of Beurling \cite{beurling70} is the following:

\begin{thm}
 A $C_0$-semigroup $(T(t))_{0\leq t<\infty}$ on a complex Banach space $\GX$ is
holomorphic  if and only if there exists a polynomial $p$ such
that 
\beq \label{eq:beurlingsuff}
\limsup_{t\to 0^+ } \|p(T(t))\|<\sup\{|p(z):|z|\leq 1\}.
\eeq
\end{thm}

Kato \cite{kato70}
and Neuberger \cite{neuberger} proved the sufficiency of (\ref{eq:beurlingsuff})
with $p(z)=z-1$, and $\sup |p(z)|=2$, providing a zero--two law for analyticity.
In general the converse is not true with $p(z)=z-1$, although it holds if $\GX$ is uniformly
convex and the semigroup is contractive \cite{pazy}.
Some extensions of this result to arbitrary Banach spaces and for semigroups which are not necessarily contractive
may be found in the very recent paper
\cite{fackler}.\\

The more recent results considered in this
survey concern estimates of the norm or spectral radius of quantities such as $T(t)-T((n+1)t)$ as $t$ tends to $0$.
These are often formulated as dichotomy results, such as the zero--quarter law (the case $n=1$ in 
the following theorem).

\begin{thm} [\cite{EM,Mok88}]\label{thm:emold}
Let $n \ge 1$ be an integer, and let $(T(t))_{t>0}$ be a semigroup in a Banach algebra. If
\[
\limsup_{t \to 0^+} \|T(t)-T((n+1)t)\| < \frac{n}{(n+1)^{1+1/n}},
\]
then either $T(t)=0$ for $t>0$ or else the closed subalgebra generated by $(T(t))_{t>0}$ is unital, and the semigroup has a bounded generator $A$: that is,
$T(t)=\exp(tA)$ for $t>0$.
\end{thm}

Another result that we mention here concerns the link between the norm and spectral radius $\rho$, and was motivated also by the
Esterle--Katznelson--Tzafriri results on estimates for $\|T^n-T^{n+1}\|$, where $T$ is a power-bounded operator (see \cite{esterle83,kt86}).
We have rewritten it in the notation of differentiable groups $(T(t)) =(\exp(tA)) $, which may even be defined for $t \in \CC$ if
$A$ is bounded; note that $T'(t)=AT(t)$.

 \begin{thm}[\cite{kalton}]\label{thm:kalton}
Let $A$ be a bounded operator on a Banach space, and let $(T(t))$ be the group given by $T(t)=\exp(tA)$. Then each of the following conditions implies that $\rho(A)=\|A\|$. \\
(i) $\sup_{t>0} t\|T'(t)\| \le 1/e$;\\
(ii) $\sup_{t>0} \|T(t)-T((s+1)t)\| \le s(s+1)^{-(1+1/s)}$ for some $s>0$;\\
(iii) $\sup_{t>0} \|T((s+i)t)-T((s-i)t)\| \le 2e^{-s \arctan(1/s)}/\sqrt{1+s^2}$ for some $s \ge 0$.
\end{thm}

The third condition of Theorem~\ref{thm:kalton} is linked to the Bonsall--Crabb proof of Sinclair's spectral radius formula for Hermitian elements of a
Banach algebra, given in \cite{bonsallcrabb}.\\

In Section~\ref{sec:23}  we review the existing literature on dichotomy laws for semigroups, first for semigroups on $\RR_+$ and then for
analytic semigroups defined on a sector in the complex plane; our main sources here are \cite{BCEP} and \cite{cep1}. Then in Section~\ref{sec:4}
we present some very recent generalizations of these results, formulated in the language of functional calculus: this discussion is based on
\cite{cep3} and \cite{cep2}.

\section{Dichotomy laws}\label{sec:23}
\subsection{Semigroups on $\RR_+$}\label{subsec:2}

We begin with a result on quasinilpotent semigroups, that is,
semigroups whose elements all have spectral radius $0$. 
This is a commonly-occurring case,  examples being 
found in the
convolution algebra $L^1(0,1)$.

\begin{thm}[\cite{Esterle}]
 Let $(T(t))_{t>0}$ be a  $C_0$-semigroup of bounded  quasinilpotent 
 linear operators on a Banach space $\GX$.  
 Then there exists $\delta>0$ such that
 $$\|T(t)-T(s)\|>\theta(s,t)\quad \mbox{ for }\quad 0<t<s<\delta,$$
 where 
 $$\theta(s,t)=(s-t)t^{t/(s-t)}s^{s/(t-s)}.$$

In particular, for all $\gamma>0$, there exists $\delta>0$ such that
 $$\|T(t)-T((\gamma+1)t)\|>\frac{\gamma}{(\gamma +1)^{1+1/\gamma}}$$
 for all $0<t<\delta$. 
\end{thm}
 
This is a sharp result, in the sense that
given a non-decreasing function
$\epsilon:(0,1)\to (0,\infty)$  there exists a nontrivial quasinilpotent
 semigroup $(T_\epsilon(t))_{t>0}$ on a Hilbert space such that:
 \[\|T_\epsilon (t)-T_\epsilon (s)\|\leq \theta(s,t)+(s-t)\epsilon(s)\]  
(see \cite{Esterle}).
Note that $\theta(s,t)=\max_{0\leq x\leq 1}(x^t-x^s)$.

It is a quantitative formulation of an intuitive fact: $T(t)$ cannot be uniformly too close to $T(s)$ for $s\neq t$, with $s,t$ small when the generator is unbounded.\\

In the non-quasinilpotent case, it is possible to
formulate similar results using the spectral radius.
The following theorem is a strengthening of
\cite[Thm~2.3]{EM}, which
is expressed in terms of
$\limsup_{t \to 0}\rho(T(t) -T(t(\g +1)))$ . Note that  $\rad\GA$ denotes the radical of the algebra $\GA$, i.e., the set of elements with   spectral radius zero.

\begin{thm}[\cite{BCEP}]\label{thm:z21}    Let $(T(t))_{t>0}$ be a non-quasinilpotent semigroup in a Banach algebra, let  $\GA$  be the closed subalgebra 
generated by $(T(t))_{t>0}$, and let $ \g>0$ be a real number. If 
there exists $t_0>0$ such that
$$\rho(T(t) -T(t(\g +1))) < \dfrac{\g}{(\g +1)^{1+\frac{1}{\g}}}$$ for $ 0<t\leq t_{0}$, then $ \GA/\rad(\GA) $ is unital, and there exist an idempotent $J$ in $\GA,$ an element $u$ of $J\GA$    and a mapping
$r : \RR_{+} \to \rad(J\GA)$, with the following properties:\\

(i)  $\phi(J) = 1$ for all $\phi \in \widehat{\GA}$;

(ii) $r(s+t) = r(s) +r(t) $ for all $s, t \in \RR_{+}$;

(iii) $JT(t) = e^{tu +r(t)}$ for $t \in \RR_{+}$, where $e^{v} = J +  \displaystyle{\sum_{k\geq1} \dfrac{v^{k}}{k!}}$ for $v \in J\GA$;

(iv) $(T(t) - JT(t))_{t\in \RR_{+}}$ is a quasinilpotent semigroup.
\end{thm}

If $\GA$ is semi-simple (that is, $\rad(\GA)=\{0\}$), then the conclusion is much more
straightforward.

\begin{cor}[\cite{BCEP}]
\label{cor:zohra}
 Let $(T(t))_{t>0}$ be a non-trivial semigroup
in a commutative semi-simple Banach algebra, let $\GA$  be the closed subalgebra generated by
$(T(t))_{t\in \RR_{+}}$ and let $ \g >0$.  If 
there exists $t_0>0$ such that
$$\rho(T(t)-T((\g + 1)t) ) < \dfrac{\g}{(\g+1)^{1+\frac{1}{\g}}}$$ 
for $ 0<t\leq t_{0}$,
then $\GA$ is unital and there exists an  element $u \in \GA$ such that $ T(t) = e^{tu}$ for $ t \in\RR_{+}$.
\end{cor}

The following theorem needs no hypothesis on $\GA$, but requires stronger estimates,
based on the norm rather than the spectral radius.

\begin{thm}[\cite{BCEP}]
\label{thm:zbigthm}
Let $(T(t))_{t>0}$ be a non-trivial semigroup in a Banach algebra, let $\GA$ be the closed
subalgebra generated by  $(T(t))_{t>0}$ and let $ n\geq1$ be an integer.
If there exists $ t_{0}>0$ such that $$\Vert(T(t) -T(t(n + 1)))\Vert < \dfrac{n}{(n +1)^{1+\frac{1}{n}}}$$ for $ 0<t \leq t_{0},$ then $\GA$ 
possesses a unit $J,\;\; \displaystyle{\lim_{t\rightarrow0^{+}}T(t)}= J$ and there exists $u \in \GA$ such that 
$T(t) = e^{tu}$ for all $t >0.$

If $(T(t)_{t>0}$ is a quasinilpotent semigroup, then 
the condition
\[\|T(t)-T((n+1)t)\|<\frac{n}{(n+1)^{1+1/n}}\quad \mbox{ for }0<t\leq t_0\]
implies that $T(t)=0$ for all $t>0$. 
\end{thm}

The sharpness of the above result is shown by the following example
\cite{BCEP}, which involves a construction of
appropriate
sequences in the non-unital Banach algebra $c_0$.

\begin{ex}
Let $G$ be an additive measurable subgroup of $\RR$ with $G\neq \RR$. Then, given 
$(\gamma_n)_n$ in $\RR^+$ such that  $t\gamma_n\in G$ for all $t\in G$
with $ t>0$ and for all $n \in \NN$, there exists 
a nontrivial semigroup
$(S(t))_{t\in G,t>0}$ in $c_0$ such that
\[ \|S(t)-S(t(\gamma_n+1))\|<\frac{\gamma_n}{(\gamma_n+1)^{1+1/\gamma_n}},      \]
for all $t\in G$, $t>0$. 

\end{ex}

\subsection{Sectorial semigroups}\label{subsec:3}

In this subsection we discuss the behaviour of analytic semigroups defined on a sector
\[
S_\alpha=\{z\in\C:\re (z)>0\mbox{ and
} |\arg(z)|<\alpha\}.
\]
with $0<\alpha \le \pi/2$. We begin with the case $\alpha=\pi/2$ and $S_\alpha=\CC_+$.

   \begin{thm}[\cite{BCEP}]\label{th:unite}
 Let $(T(t))_{t\in \C_+}$ be an analytic non-quasinilpotent semigroup in a Banach algebra. Let ${\cal A}$ be the closed subalgebra  generated by $(T(t))_{t\in \C_+}$ and let $\gamma >0.$
   If there exists $t_0>0$ such that 
   
   $$\sup_{t \in \C_+, \vert t \vert \le t_0}\rho(T(t)-T(\gamma+1)t))<2$$
   then ${\cal A}/\rad {\cal A}$ is unital, and the generator of $(\pi(T(t)))_{t>0}$ is bounded, where $\pi:{\cal A}\to {\cal A}/\rad {\cal A}$ denotes the canonical surjection.
   \end{thm}

  A semigroup $(T(t))$ defined on the positive reals or on a sector is said to be {\em exponentially bounded\/} if there exist $c_1>0$ and $c_2\in \R$ such that $\Vert T(t)\Vert \le c_1e^{c_2\vert t \vert}$ for every $t.$   Beurling~\cite{beurling70}, in his work described in the introduction,
showed that there exists a universal constant $k$ such that every exponentially bounded weakly measurable semigroup $(T(t))_{t>0}$ of bounded operators satisfying 
\[
\limsup_{t\to 0^+}\Vert \Id -T(t)\Vert =\rho <2
\]
 admits an exponentially bounded analytic extension to a sector $S_{\alpha}$ with $\alpha \ge k(2-\rho)^2$.
From this one easily obtains the following result.

\begin{thm}[\cite{BCEP}] Let $(T(t))_{t\in \C_+}$ be an analytic semigroup of bounded operators on a Banach space $\GX.$ If the generator of the semigroup is unbounded, then we have, for $-{\pi \over 2} <\alpha < {\pi\over 2}$,
\[
\limsup_{t \to 0^+}\Vert \Id -T(t)\Vert \ge 2 -\sqrt{{\pi\over 2}- \vert \alpha \vert\over k},
\]
where $k$ is Beurling's universal constant.
\end{thm} 

We now consider similar results on smaller sectors than the
half-plane, and in fact the result we prove will be stated in a far
more general context.

\begin{thm}[\cite{BCEP}]\label{thm:moregeneral}
Let  $0 < \alpha < \pi/2$ and let $f$ be an entire function with $f(0)=0$ and $f(\RR) \subseteq \RR$, such that
\begin{equation}\label{eq:vanish}
 \sup_{\re z > r} |f(z)| \to 0 \qquad \hbox{as} \quad r \to \infty,
\end{equation}
and $f$ is a linear combination of functions of the form $z^m \exp(-zw)$ for $m =0,1,2,\ldots$
and $w >0$.
Let $(T(t))_{t\in S_\alpha}=(\exp(tA))_{t \in S_\alpha}$ be an analytic non-quasinilpotent semigroup in
a Banach algebra and let  
 $\cal A$ be the subalgebra generated by
$(T(t))_{t\in S_\alpha}$.
If there exists $t_0>0$ such that 
\[ \sup_{t\in S_\alpha,|t|\le t_0}\rho (f(-tA))<k(S_\alpha),\]
with $k(S_\alpha)=\sup_{t\in S_\alpha}|f(z)|$,
then ${\cal A}/\rad{\cal A}$ is unital
and the generator of $\pi(T(t))_{t\in S_{\alpha}}$ is bounded, where $\pi: {\cal A}\to {{\cal A}/\rad({\cal A})}$ denotes the canonical surjection.
\end{thm}

Note that $f(-tA)$ is well-defined in terms of $T(t)$ and its derivatives.\\

 Suitable examples of $f(z)$ are linear combinations of
 $z^m \exp(-z)$, $m=1,2,3,\ldots$, and
$\exp(-z)-\exp(-(\gamma+1)z)$; also real linear combinations of the form
$\sum_{k=1}^n a_k \exp(-b_k z)$ with $b_k>0$ and $\sum_{k=1}^n a_k=0$.
This provides results analogous to those of \cite[Thm.~4.12]{kalton},
where the behaviour of expressions such as $\|tA\exp(tA)\|$ and $\|\exp(tA)-\exp(stA)\|$ was considered
for all $t>0$.

\begin{rem}
Another function considered in \cite{kalton} is $f(z)=e^{-sz}\sin z$, where we now require
$s > \tan\alpha$ for $f(-tA)$ to be well-defined for $t \in S_\alpha$. This does not satisfy the condition
(\ref{eq:vanish}), but we note that it holds for $z \in S_\alpha$, while for
$z \not\in S_\alpha$ there exists a constant $C>0$ 
with the following property: 
for each $z$ with $\re z>C$ there exists $\lambda\in(0,1)$ such that
$|f(\lambda z)|\ge \sup_{z \in S_\alpha} |f(z)|$.
Using this observation, it is  not difficult 
to adapt the proof of Theorem \ref{thm:moregeneral} to this case.
\end{rem}

The sharpness of the constants can be shown by considering examples in $C_0([0,1])$.

One particular case of the above is used in the estimates considered
by Bendaoud, Esterle and 
Mokhtari \cite{BEM,EM}.

\begin{cor}\label{cor:sector2}
Let $\gamma>0$ and $0 < \alpha < \pi/2$.
Let $(T(t))_{t\in S_\alpha}$ be an analytic non-quasinilpotent semigroup in
a Banach algebra and let  
 $\cal A$ be the closed subalgebra generated by
$(T(t))_{t\in S_\alpha}$.
If there exists $t_0>0$ such that 
\[ \sup_{t\in S_\alpha,|t|\le t_0}\rho (T(t)-T(t(\gamma
+1)))<k(S_\alpha),\]
with $k(S_\alpha)=\sup_{t\in S_\alpha}|\exp(-t)-\exp(-(\gamma+1)t)|$,
then ${\cal A}/\rad{\cal A}$ is unital
and the generator of $\pi(T(t))_{t\in S_{\alpha}}$ is bounded,
where $\pi: {\cal A}\to {{\cal A}/\rad({\cal A})}$ denotes the canonical surjection.  
\end{cor}

Now set $f_n(z)=z^ne^{-z},$ and set $k_n(\alpha)=\max_{z \in S_{\alpha}}|f_n(z)|$. A
straightforward  computation shows that $k_n(\alpha)=\left ( {n \over e \cos(\alpha)}\right )^n.$

If $A$ is the generator of an analytic semigroup $(T(t))_{t\in S_{\alpha}}$, 
then we have $f_n(-tA)=(-1)^n t^nT^{(n)}(t).$ So the following result,
which  may be deduced from Hille's work,
described in Theorem~\ref{thm:hille50},
means that if 
\[
\sup_{t\in S_{\alpha}, 0< \vert t \vert <\delta}\Vert f_n(-tA)\Vert <k_n(\alpha)
\]
 for
some $\delta >0,$ then the generator of the semigroup is bounded.

\begin{thm}[\cite{BCEP}] Let $n \ge 1$ be an integer, let $\alpha \in (0,\pi/2)$ and let $(T(t))_{t\in S_{\alpha}}$ be an analytic semigroup. If \[
 \sup_{t\in S_{\alpha}, 0< \vert t \vert <\delta}\Vert t^nT^{(n)}(t)\Vert<\left ( {n \over e \cos(\alpha)}\right )^n
\]
 for some $\delta >0,$ then the closed algebra generated by the semigroup is unital, and the generator of the semigroup is bounded.
\end{thm}

The remainder of this section is devoted to quasinilpotent semigroups.
We let
$D(0,r)$ denote $\{z \in \CC: |z|<r\}$.

\begin{rem}\label{rem:es1}
An analytic semigroup $(T(t))_{t \in S_\alpha}$ 
acting on a Banach space $\GX$ and
bounded near the origin can
be extended to the closed sector $\overline{S_\alpha}$. 
Indeed, assume that there exists $r >0$ such that 
\[\sup_{t \in D(0,r)\cap S_\alpha}\Vert T(t) \Vert <+\infty.\]
Then  $\lim_{\stackrel {t \to w}{_{t \in S_\alpha}}}T(t)x$ exists for every 
$x \in \GX$ and every $w \in \partial S_\alpha$. Moreover if we set
\[
T(w)x=\lim_{\stackrel {t \to w}{_{t \in S_\alpha}}}T(t)x,
\]
then $(T(t))_{t \in \overline{S_\alpha}}$ is a semigroup of bounded operators which
is continuous with respect to the strong operator topology.
For we have $\lim_{\stackrel {t \to w}{_{t \in S_\alpha}}}T(t)T(t_0)x=T(t_0)x$ for every $t_0 >0$ and every $x \in \GX$. Now the result follows immediately from the
 fact that $\bigcup_{t>0}T(t)\GX$ is dense in $\GX$, given that 
\[
\sup_{z \in D(0,r)\cap S_\alpha}\Vert T(t) \Vert <+\infty.
\] 
\end{rem}

The next lemma demonstrates that  nontrivial 
quasinilpotent analytic semigroups cannot be bounded on the right 
half-plane $\CC_+$. In fact, more is true. 

\begin{lem}[\cite{cep1}]\label{lem:esbounded}
Let $(T(t))_{t \in \CC_+}$ be a quasinilpotent analytic semigroup of bounded operators on a Banach space $\GX$.
Suppose that there exists $r >0$ such that 
\[
\sup_{t \in D(0,r)\cap \CC_+}\Vert T(t) \Vert <+\infty,
\] 
and
define $T(iy)$ for $y \in \RR$ using Remark~\ref{rem:es1}.
 If 
\[\int_{-\infty}^{\infty} \frac{\log^+\|T(iy)\|}{1+y^2} \, dy <+\infty,\]
 then $T(t)=0$ for $t \in \CC_+$.
\end{lem}

In the case when the semigroup is bounded near the origin, we may give appropriate estimates on the imaginary axis.

\begin{thm}[\cite{cep1}]\label{th:esdist}
  Let $(T(t))_{t \in \CC_+}$ be a nontrivial quasinilpotent analytic semigroup satisfying the conditions of Remark~\ref{rem:es1}, 
and let $s>0$. Then 
\[\max(\rho(T(iy)-T(iy+is)), \rho(T(-iy)-T(-iy-is)))\ge 2,\]
for every $y >0$. 
\end{thm}

From this we may obtain estimates for semigroups satisfying
a growth condition near the imaginary axis.

\begin{cor} [\cite{cep1}] Let $(T(t))_{t \in \CC_+}$ be a quasinilpotent analytic semigroup such that
\[
\sup_{y \in \RR} e^{-\mu \vert y \vert}\Vert T(\delta +iy)\Vert < +\infty
\] for some $\delta >0$ and some $\mu >0$, and let $\gamma >0$.  Then 
\[\sup_{t \in D(0,r)\cap \CC_+}\Vert T(t)-T((1+\gamma) t)\Vert\ge 2,\]
for every $r>0.$
\end{cor}

\section{Lower estimates for functional calculus}\label{sec:4}

In this section we summarise some very recent results from \cite{cep3} 
and \cite{cep2}, which provide  far-reaching generalizations of earlier work.

\subsection{Semigroups on $\RR_+$}

\subsubsection{The quasinilpotent case}

Recall that if $(T(t))_{t>0}$ is a  uniformly bounded  strongly continuous semigroup with generator $A$, then
\beq\label{eq:41rev}
(A+\lambda \Id)^{-1}= - \int_0^\infty e^{\lambda t} T(t) \, dt,
\eeq
for all $\lambda \in \CC$ with $\re \lambda < 0$.
Here the integral is taken in the sense of Bochner with respect to the strong operator topology.
If, in addition, $(T(t))_{t > 0}$ is quasinilpotent, then
we have (\ref{eq:41rev})
for all $\lambda \in \CC$.

Similarly, if $\mu \in M_c(0,\infty)$ (the space of complex finite
Borel measures
on $(0,\infty)$) with Laplace transform
\beq\label{eq:Lmu}
F(s):=\LL \mu(s) = \int_0^\infty e^{-s\xi} \, d\mu(\xi) \qquad (s \in \CC_+),
\eeq
and $(T(t))_{t>0}$ is a strongly continuous semigroup of bounded operators on 
a Banach space $\GX$, then we have a functional
calculus for its generator $A$, defined by
\[
F(-A)= \int_0^\infty T(\xi) \, d\mu(\xi),
\]
in the sense of the strong operator topology; i.e.,
\[
F(-A)x= \int_0^\infty T(\xi)x \, d\mu(\xi)  \qquad (x \in \GX),
\]
which exists as a Bochner integral.

The following theorem applies to several examples studied  in \cite{BCEP,Esterle,EM,kalton}; these include $\mu=\delta_1-\delta_2$, the difference of two Dirac
measures, where $F(s):=\LL \mu(s)=e^{-s}-e^{-2s}$ and $F(-sA)=T(t)-T(2t)$. 
More importantly, the theorem applies to many other examples, such as $d\mu(t)=(\chi_{[1,2]}-\chi_{[2,3]})(t) dt$ and $\mu=\delta_1-3\delta_2+\delta_3+\delta_4$,
which are not accessible with the methods of \cite{BCEP,Esterle,EM,kalton}.

\begin{thm}\label{thm:thmhelp}
Let $\mu \in M_c(0,\infty)$ be a real measure such that $\ds \int_0^\infty d\mu(t) = 0$, and let
$(T(t))_{t > 0}$ be a strongly continuous quasinilpotent semigroup of bounded operators on a Banach space $\GX$. Set $F=\LL \mu$. Then
there exists $\eta > 0$ such that
\[
\| F(-sA ) \|  > \max_{x \ge 0} |F(x)| \qquad \hbox{for} \quad 0 < s \le \eta.
\]
\end{thm}

If $\mu \in M_c(0,\infty)$ is a complex measure, then we write $\widetilde F=\LL \overline\mu$, so that
$\widetilde F(z)=\overline{F(\bar z)}$. By considering the real measure $\nu:=\mu * \overline\mu$,
we obtain the following result.

\begin{cor}
Let $\mu \in M_c(0,\infty)$ be a complex measure such that $\ds \int_0^\infty d\mu(t) = 0$, and let
$(T(t))_{t > 0}$ be a strongly continuous quasinilpotent semigroup of bounded operators on a Banach space $\GX$. Set $F=\LL \mu$. Then
there exists $\eta > 0$ such that
\[
\| F(-sA )\widetilde F(-sA) \|  > \max_{x \ge 0} |F(x)|^2 \qquad \hbox{for} \quad 0 < s \le \eta.
\]
\end{cor}

\subsubsection{The non-quasinilpotent case}

Recall that 
a sequence $(P_n)_{n \ge 1}$
of idempotents in a Banach algebra $\GA$ is said to   be {\em exhaustive \/} 
if $P_n^2=P_nP_{n+1} = P_n$ for all $n$ and if
for every $\chi \in \widehat \GA$ there is a $p$ such that
$\chi(P_n)=1$ for all $n \ge p$. Such sequences may often be
found in non-unital algebras: for example, $P_n=e_1+\ldots+e_n$
($n=1,2,\ldots$)
in the Banach algebra $c_0$.

\begin{thm}\label{thm:dodgy}
Let $(T(t))_{t>0}$ be a strongly continuous and eventually norm-continuous non-quasinilpotent semigroup on a Banach space $\GX$, with generator $A$. Let 
$F=\LL \mu$, where $\mu \in M_c(0,\infty)$ is a real measure such that $\int_0^\infty d\mu=0$. If there exists
$(u_k)_k \subset (0,\infty)$ with $u_k \to 0$ such that
\[
\rho(F(-u_k A)) < \sup_{x >0} |F(x)|,
\]
then the algebra $ \GA$ 
 generated by
$(T(t))_{t>0}$
possesses an exhaustive sequence of idempotents $(P_n)_{n \ge 1}$
such that each semigroup $(P_n T(t))_{t>0}$ has a bounded generator.

If, further, $\|F(-u_k A)\| < \sup_{x>0} |F(x)|$, then
$\bigcup_{n \ge 1} P_n \GA$ is dense in $\GA$.
\end{thm}

\subsection{Analytic semigroups}

For $0<\alpha< \pi/2$, let $H(S_\alpha)$ denote the Fr\'echet space of holomorphic functions on $S_\alpha$, endowed with the topology of local uniform convergence; thus, if $(K_n)_{n \ge 1}$ is an increasing sequence of
compact subsets of $S_\alpha$ with $\bigcup_{n \ge 1} K_n=S_\alpha$, we may specify the topology by the seminorms
\[
\|F\|_n:= \sup \{|f(z)|: z  \in K_n\}.
\]
We write $H(S_\alpha)'$ for its dual space; that is, the space of
continuous linear functionals
$\phi: H(S_\alpha) \to \CC$.
This means that there is an index $n$ and a constant $M>0$ such that 
$|\langle f,\phi \rangle| \le M \|f\|_n$ for all $f \in H(S_\alpha)$. 

We define the {\em Fourier--Borel transform\/} of $\phi$ by
\[
\FB(\phi)(z) = \langle e_{-z}, \phi \rangle,
\]
for $z \in \CC$, where $e_{-z}(\xi)=e^{-z\xi}$ for $\xi \in S_\alpha$.

If $\phi \in H(S_\alpha)'$, as above, then by the Hahn--Banach theorem, it can be extended to a functional on $C(K_n)$, which we still write as $\phi$, and is thus
induced by a Borel measure $\mu$ supported on $K_n$.

That is, we have
\[
\langle f,\phi\rangle = \int_{S_\alpha} f(\xi) \, d\mu(\xi),
\]
where $\mu$ (which is not unique) is a compactly supported measure. For example, if $\langle f,\phi \rangle = f'(1)$, then
\[
\langle f,\phi \rangle = \frac{1}{2\pi i} \int_C \frac{f(z) \, dz}{(z-1)^2},
\]
where $C$ is any sufficiently small circle surrounding the point 1.
Note that 
\[
\FB(\phi)(z)=\int_{S_\alpha} e^{-z\xi} d\mu(\xi).
\]

Now let $T:=(T(t))_{t \in S_\alpha}$ be an analytic semigroup  on a  Banach space $\GX$, with  generator $A$. Let $\phi\in H(S_\alpha)'$ and let $F=\FB(\phi)$.

We may thus define, formally to start with,
\[
F(-A)=\langle T,\phi \rangle = \int_{S_\alpha} T(\xi) \, d\mu(\xi),
\]
which is well-defined as a Bochner integral in $\GA$.
It is easy to verify that the definition is independent of the choice of $\mu$ representing $\phi$.

Indeed, if $u \in S_{\alpha-\beta}$, where $\supp \mu \subset S_\beta$ and $0<\beta<\alpha$, then
we may also define 
\[
F(-uA) = \int_{S_\beta} T(u\xi) \, d\mu(\xi),
\]
since $u\xi$ lies in $S_\alpha$.

The following theorem extends \cite[Thm. 3.6]{BCEP}. In the following, a symmetric measure is a measure such that
$\mu(\overline S)=\overline{\mu(S)}$ for $S \subset S_\alpha$.  A symmetric measure will have a Fourier--Borel transform $f$ 
satisfying $f(z)=\widetilde f(z):=\overline{f(\overline z)}$ for all $z \in \CC$.

\begin{thm}\label{thm:moregeneral2}
Let  $0 < \alpha < \pi/2$.
Let $\phi \in H(S_\alpha)'$, induced by a symmetric measure $\mu \in M_c(S_\alpha)$ such that $\int_{S_\alpha} d\mu(z)=0$,
and let $f=\FB( \phi)$.  
Let $(T(t))_{t\in S_\alpha}=(\exp(tA))_{t \in S_\alpha}$ be an analytic non-quasinilpotent semigroup in
a Banach algebra and let  
 $\GA$ be the subalgebra generated by
$(T(t))_{t\in S_\alpha}$.
If there exists $t_0>0$ such that 
\[ \sup_{t\in S_\alpha,|t|\le t_0}\rho (f(-tA))< \sup_{z\in S_\alpha}|f(z)|,\]
then ${\GA}/\rad{\GA}$ is unital
and the generator of $\pi(T(t))_{t\in S_{\alpha}}$ is bounded, where $\pi: {\GA}\to {{\GA}/\rad({\GA})}$ denotes the canonical surjection.
\end{thm}

By considering the convolution of a functional $\phi \in H(S_\alpha)'$, with Fourier--Borel transform $f$, and 
the functional $\widetilde \phi$ with Fourier--Borel transform $\widetilde f$, we obtain the following result.

\begin{cor}
Let  $0 < \alpha < \pi/2$.
Let $\phi \in H(S_\alpha)'$, induced by a  measure $\mu \in M_c(S_\alpha)$ such that $\int_{S_\alpha} d\mu(z)=0$,
and let $f=\FB( \phi)$.
Let $(T(t))_{t\in S_\alpha}=(\exp(tA))_{t \in S_\alpha}$ be an analytic non-quasinilpotent semigroup in
a Banach algebra and let  
 $\GA$ be the subalgebra generated by
$(T(t))_{t\in S_\alpha}$.
If there exists $t_0>0$ such that 
\[ \sup_{t\in S_\alpha,|t|\le t_0}\rho (f(-tA)\widetilde f(- t A))< \sup_{z\in S_\alpha}|f(z)||\widetilde f(z)|,\]
then ${\GA}/\rad{\GA}$ is unital
and the generator of $\pi(T(t))_{t\in S_{\alpha}}$ is bounded, where $\pi: {\GA}\to {{\GA}/\rad({\GA})}$ denotes the canonical surjection.
\end{cor}

It is possible to obtain a similar conclusion, based only on estimates on the positive real line.

\begin{thm}\label{thm:4.6}
Let  $0 < \alpha < \pi/2$.
Let $\phi \in H(S_\alpha)'$, induced by a  symmetric measure $\mu \in M_c(S_\alpha)$ such that $\int_{S_\alpha} d\mu(z)=0$,
and let $f=\FB( \phi)$. Suppose that 
$f(\RR_+) \subset \RR$.
Let $(T(t))_{t\in S_\alpha}=(\exp(tA))_{t \in S_\alpha}$ be an analytic non-quasinilpotent semigroup in
a Banach algebra and let  
 $\GA$ be the subalgebra generated by
$(T(t))_{t\in S_\alpha}$.
If there exists $t_0>0$ such that 
\[ \rho (f(-tA))< \sup_{x>0}|f(x)|,\]
for all $0 < t \le t_0$,
then ${\GA}/\rad{\GA}$ is unital
and the generator of $\pi(T(t))_{t\in S_{\alpha}}$ is bounded, where $\pi: {\GA}\to {{\GA}/\rad({\GA})}$ denotes the canonical surjection.
\end{thm}

The following example shows that the hypotheses of Theorem~\ref{thm:4.6} are sharp. 

\begin{ex} In the Banach algebra
$\GA=C_0([0,1])$ consider the semigroup $T(t):x \mapsto x^t$. Clearly, 
$(T(t))$ is not norm-continuous at $0$.

For $x \in (0,1]$ (which can be identified with the Gelfand space of $\GA$) let $f=\FB(\mu)$ and
\[
f(-tA)(x) = \int_{S_\alpha} x^{-t\xi} \, d\mu(\xi) = \int_{S_\alpha} e^{-t\xi \log x} \, d\mu(\xi),
\]
where $\mu \in M_c(S_\alpha)$, supposing that  $\int_{S_\alpha} d\mu(z)=0$ and that $f(\RR_+) \subset \RR$.

Thus $f(-tA)(x)=f(-t \log x)$ and 
\[
\rho(f(-tA))=\|f(-tA)\|= \sup_{x>0}|f(-t\log x)|=\sup_{r>0}|f(tr)|.
\]
Clearly,  
\[ \sup_{t\in S_\alpha,|t|\le t_0}\rho (f(-tA))= \sup_{t\in S_\alpha}|f(z)|
\]
for all $t_0>0$.
\end{ex}

\section*{Acknowledgements}

This work was partially supported by the 
London Mathematical Society (Scheme 2).


\begin{thebibliography}{10}

\bibitem{BCEP}
Z.~Bendaoud, I.~Chalendar, J.~Esterle, and J.~R. Partington.
\newblock Distances between elements of a semigroup and estimates for
  derivatives.
\newblock {\em Acta Math. Sin. (Engl. Ser.)}, 26(12):2239--2254, 2010.

\bibitem{BEM}
Z.~Bendaoud, J.~Esterle, and A.~Mokhtari.
\newblock Distances entre exponentielles et puissances d'\'el\'ements de
  certaines alg\`ebres de {B}anach.
\newblock {\em Arch. Math. (Basel)}, 89(3):243--253, 2007.

\bibitem{beurling70}
A.~Beurling.
\newblock On analytic extension of semigroups of operators.
\newblock {\em J. Functional Analysis}, 6:387--400, 1970.

\bibitem{bonsallcrabb}
F.~F. Bonsall and M.~J. Crabb.
\newblock The spectral radius of a {H}ermitian element of a {B}anach algebra.
\newblock {\em Bull. London Math. Soc.}, 2:178--180, 1970.

\bibitem{cep1}
I.~Chalendar, J.~Esterle, and J.~R. Partington.
\newblock Boundary values of analytic semigroups and associated norm estimates.
\newblock In {\em Banach algebras 2009}, volume~91 of {\em Banach Center
  Publ.}, pages 87--103. Polish Acad. Sci. Inst. Math., Warsaw, 2010.

\bibitem{cep3}
I.~Chalendar, J.~Esterle, and J.~R. Partington.
\newblock In preparation, 2014.

\bibitem{cep2}
I.~Chalendar, J.~Esterle, and J.~R. Partington.
\newblock Lower estimates near the origin for functional calculus on operator
  semigroups.
\newblock Submitted, 2014.

\bibitem{esterle83}
J.~Esterle.
\newblock Quasimultipliers, representations of {$H\sp{\infty }$}, and the
  closed ideal problem for commutative {B}anach algebras.
\newblock In {\em Radical {B}anach algebras and automatic continuity ({L}ong
  {B}each, {C}alif., 1981)}, volume 975 of {\em Lecture Notes in Math.}, pages
  66--162. Springer, Berlin, 1983.

\bibitem{Esterle}
J.~Esterle.
\newblock Distance near the origin between elements of a strongly continuous
  semigroup.
\newblock {\em Ark. Mat.}, 43(2):365--382, 2005.

\bibitem{EM}
J.~Esterle and A.~Mokhtari.
\newblock Distance entre \'el\'ements d'un semi-groupe dans une alg\`ebre de
  {B}anach.
\newblock {\em J. Funct. Anal.}, 195(1):167--189, 2002.

\bibitem{fackler}
S.~Fackler.
\newblock Regularity of semigroups via the asymptotic behaviour at zero.
\newblock {\em Semigroup Forum}, 87(1):1--17, 2013.

\bibitem{hille50}
E.~Hille.
\newblock On the differentiability of semi-group operators.
\newblock {\em Acta Sci. Math. Szeged}, 12:19--24, 1950.

\bibitem{HP74}
E.~Hille and R.~S. Phillips.
\newblock {\em Functional analysis and semi-groups}.
\newblock American Mathematical Society, Providence, R. I., 1974.
\newblock Third printing of the revised edition of 1957, American Mathematical
  Society Colloquium Publications, Vol. XXXI.

\bibitem{kalton}
N.~Kalton, S.~Montgomery-Smith, K.~Oleszkiewicz, and Y.~Tomilov.
\newblock Power-bounded operators and related norm estimates.
\newblock {\em J. London Math. Soc. (2)}, 70(2):463--478, 2004.

\bibitem{kato70}
T.~Kato.
\newblock A characterization of holomorphic semigroups.
\newblock {\em Proc. Amer. Math. Soc.}, 25:495--498, 1970.

\bibitem{kt86}
Y.~Katznelson and L.~Tzafriri.
\newblock On power bounded operators.
\newblock {\em J. Funct. Anal.}, 68(3):313--328, 1986.

\bibitem{Mok88}
A.~Mokhtari.
\newblock Distance entre \'el\'ements d'un semi-groupe continu dans une
  alg\`ebre de {B}anach.
\newblock {\em J. Operator Theory}, 20(2):375--380, 1988.

\bibitem{neuberger}
J.~W. Neuberger.
\newblock Analyticity and quasi-analyticity for one-parameter semigroups.
\newblock {\em Proc. Amer. Math. Soc.}, 25:488--494, 1970.

\bibitem{pazy}
A.~Pazy.
\newblock {\em Semigroups of linear operators and applications to partial
  differential equations}, volume~44 of {\em Applied Mathematical Sciences}.
\newblock Springer-Verlag, New York, 1983.

\end{thebibliography}

\end{document}